\documentclass[12pt]{article}

\usepackage{amsmath}
\usepackage{multirow} 
\usepackage{longtable}
\usepackage{latexsym,amsfonts}
\usepackage{graphicx}
\usepackage{subcaption}
\usepackage{pdflscape}
\usepackage{caption}
\usepackage{tikz-cd} 
\usepackage{tikz}
\usepackage{enumitem} 
\usepackage{authblk}
\newtheorem{thm}{Theorem}[section]

\newtheorem{lem}[thm]{Lemma}

\author[1]{Sanja Rukavina}
\author[2]{Vladimir D. Tonchev}
\affil[1]{Faculty of Mathematics, University of Rijeka, 51000 Rijeka, Croatia\\

E-mail: sanjar@math.uniri.hr, ORCID: 0000-0003-3365-7925}
\affil[2]{Department of Mathematical Sciences, Michigan Technological University, Houghton, MI 49931, USA\\

E-mail: tonchev@mtu.edu, ORCID: 0000-0003-1806-3571
}
\title{Symmetric $2$-$(35,17,8)$ designs with an automorphism of order 2}

\begin{document} 

\maketitle

\begin{abstract}

The largest prime $p$ that can be the order of an automorphism 
of a 2-$(35,17,8)$ design is $p=17$, and all 2-$(35,17,8)$ designs with an automorphism of order 17 
were classified in  \cite{Ton86}. 
The symmetric 2-$(35,17,8)$
designs with automorphisms of an odd prime order $p<17$ were classified in  \cite{crn_r35} and  \cite{BFW}. 
In this paper we give the classification of all symmetric 
$2$-$(35,17,8)$ designs that admit an automorphism of order $p=2$.  It is shown that there are exactly $11,642,495$ 
nonisomorphic such designs.  
Furthermore, it is shown that the number of nonisomorphic  $3$-$(36,18,8)$ designs
which have at least one derived 2-$(35,17,8)$ design with an automorphism of order 2, is  $1,015,225$.

\end{abstract}



\section {Introduction}  

We assume familiarity with basic facts and notions from the theory of combinatorial designs \cite{BJL,CRC, l, ton88}.
 
The parameters 2-$(35,17,8)$ are the smallest parameters
for  a symmetric design such that a complete classification up to isomorphism
of all such designs is yet unknown. It was proved in \cite{Ton86} that the largest prime $p$ that can be 
the order of an automorphism  of a 2-$(35,17,8)$ design,  is $p=17$, and that there are exactly 21 nonisomorphic designs with an automorphism of order 17. 
Further, it was proved in \cite{crn_r35} that $p<17$ implies $p \in \{ 2, 3, 5, 7 \}$,  and all  2-$(35,17,8)$ designs with an automorphism of order $5$ or $7$ were classified. Finally, in \cite{BFW},  Bouyukliev, Fack and Winne completed the classification of all 2-$(35,17,8)$ designs with an automorphism
of an odd prime order $p<17$, and gave a partial enumeration of symmetric $2$-$(35,17,8)$ designs with an automorphism of order 2. The results of this partial enumeration are summarized in Table \ref{bfw_table}.

In this paper, we give a complete classification of all symmetric $2$-$(35,17,8)$ designs with an automorphism  of order 2.
The classification process consists of four steps: 1) finding the possible actions of a group of order 2 on the sets of
points and blocks, 2) finding the corresponding orbit matrices, 3)  expanding (or {\it indexing} in the terminology of 
 \cite{c-r,cep,jan})
the orbit matrices to  (0,1)-incidence matrices of designs, and 4) classifying the resulting designs up to isomorphism.
In Section 2, we describe all possible actions of an automorphism of order 2 on a symmetric $2$-$(35, 17, 8)$ design, 
and give a detailed classification of all such symmetric designs.
In Section 3, we use  these results to classify all $3$-$(36,18,8)$ designs that have at least one derived
2-$(35,17,8)$  design an automorphism of order 2.

 \begin{table}[htpb!]
\begin{center} \begin{footnotesize}
\begin{tabular}{|c ||c| c | c|c|c|c|c|c|c| }
 \hline 

Number of fixed points& 1& 3 & 5 & 7 & 9 & 11 & 13 & 15 & 17 \\
 \hline
Number of designs&0&111098&237058& ? & 430656 & ? &? &?&?   \\
\hline
  
\end{tabular} \end{footnotesize}
 \caption{A partial classification of symmetric $2$-$(35,17,8)$ designs with an automorphism of order 2  \cite{BFW}}\label{bfw_table}
\end{center} 
\end{table}

\section{Classification of symmetric $2$-$(35,17,8)$  designs with an automorphism of order 2}
\label{d35}

Let $ {\mathcal D}=( {\mathcal P},{\mathcal B}) $ be a symmetric $2$-$ (v,k,\lambda ) $ design with  set of points $\mathcal P$ and set of blocks $\mathcal B$. An automorphism of $ {\mathcal D}=( {\mathcal P},{\mathcal B}) $  is a permutation on the set $\mathcal P$ that maps blocks to blocks. 
The set of all automorphisms of ${\mathcal D}$ forms its full automorphism group, denoted by $Aut({\mathcal D})$. 
The action of an automorphism of order 2 on ${\mathcal D}$ generates the same number of point and block orbits (see \cite[Theorem 3.3]{l}). We denote this number by $t$, the point orbits by ${\cal P}_{1},\ldots,{\cal P}_{t}$,
the block orbits by $ {\cal B}_{1},\ldots ,{\cal B}_{t} $, and set $ |{\cal P}_{r}|=\omega _{r} $ and $ |{\cal B}_{i}|=\Omega _{i} $.
Let $ \gamma _{ir} $ denote the number of points of $ {\cal P}_{r} $ which are incident with a representative of the block orbit ${\cal B}_{i} $.
The following equations apply to these numbers (cf. e.g. \cite{c-r}):
\begin{eqnarray}  
\sum _{r=1}^{t}{\gamma }_{ir} & = & k \label{eq-1}\, ,\\
\sum _{r=1}^{t}\frac{{\Omega }_{j}}{{\omega }_{r}}{\gamma }_{ir}{\gamma }_{jr} & = & \lambda {\Omega }_{j}+{\delta }_{ij}\cdot (k-\lambda )\label{eq-2}\,.
\end{eqnarray}

A $(t \times t)$-matrix $({\gamma}_{ir})$ with entries satisfying the  conditions $(1)$ and $(2)$ is called an orbit matrix for
the parameters
$(v,k, \lambda)$ and orbit lengths
distributions $(\omega_{1}, \ldots ,\omega_{t})$,
$(\Omega_{1}, \ldots ,\Omega_{t})$.\\

The first step in the classification process is to determine all possible actions of an automorphism of order 2 on a symmetric $2$-$(35,17,8)$ design, and then construct the corresponding orbit matrices that satisfy the equations (\ref{eq-1}) and (\ref{eq-2}). The resulting orbit matrices are then expanded into (0,1)-incidence matrices of symmetric $2$-$(35,17,8)$ designs by the indexing
method described by Janko \cite{jan}. Finally, an isomorphism test is performed on the obtained designs. 
In our constructions, we use computers. In addition to our own computer programs, we use a computer program by  \'{C}epuli\'{c} \cite{cep} for the construction of orbit matrices, the System for Computational Discrete Algebra GAP \cite{GAP2022}, and its package DESIGN \cite{Soi23} to compute the automorphism groups and the representatives of the isomorphism classes of the constructed designs.\\

Denote by $f$ the number of fixed points for an action of an automorphism of order 2 on a $2$-$(35,17,8)$ design. Then $3 \leq f \leq 17$ (see \cite[Corollary 3.7, Proposition 4.23]{l}). The following lemma excludes the cases for which there are no orbit matrices.

\begin{lem}
Let $f$ be the number of fixed points of an automorphism of order 2 acting on a symmetric $2$-$(35,17,8)$ design. 
Then $f \notin \{13,17 \}$.
\end{lem}

{\it Proof.}  As mentioned in \cite{BFW}, the main obstacle in constructing orbit matrices for an automorphism of order two acting with a large number of fixed points is the construction of the fixed part for such orbit matrices. Therefore, for the cases where $f\geq 13$ we start our construction by considering the rows of orbit matrices that correspond to the block orbits of length 2.

Suppose that $f=13$. In this case, an orbit matrix should have $11$  rows corresponding to  block orbits of length 2, 
and there are only three possible types for these rows (the first $13$ entries in each type correspond to the fixed points,
e.g., the point orbits of length 1):
\begin{center}
1  1  1  1  1  1  1  1  0  0  0  0  0  1  1  1  1  1  1  1  1  1  0  0 \\
1  1  1  1  1  1  0  0  0  0  0  0  0  2  1  1  1  1  1  1  1  1  1  0 \\
1  1  1  1  0  0  0  0  0  0  0  0  0  2  2  1  1  1  1  1  1  1  1  1 \\
\end{center}

It is easy to check that by using rows of these types it is not possible to construct more than $7$ rows satisfying (\ref{eq-2}).  This implies the non-existence of orbit matrices.\\

Suppose now that $f=17$. In this case, there is only one type for rows of an orbit matrix corresponding to 
the block orbits of size 2:
\begin{center}
1  1  1  1  1  1  1  1  0  0  0  0  0  0  0  0  0  1  1  1  1  1  1  1  1  1 
 \end{center}

 It is easy to check that it is not possible to construct two rows of this type  that satisfy (\ref{eq-2}).

It follows that orbit matrices do not exist for an automorphism of order 2 acting on a symmetric $2$-$(35,17,8)$  
design with $13$ or $17$ fixed points. $\Box$\\

 The number of constructed orbit matrices for the other cases, which are available at
\begin{verbatim}
 https://www.math.uniri.hr/~sanjar/structures/
 \end{verbatim}
is given in   the second row of Table \ref{OM}.  The results of further analysis of these orbit matrices, which consists in constructing corresponding symmetric $2$-$(35,17,8)$  designs and elimination of isomorphic copies, are also shown in Table \ref{OM}. The third row of Table \ref{OM} contains the total number of non-isomorphic  $2$-$(35,17,8)$  designs for given action of an automorphism of order 2. The fourth row gives information about how many of the constructed symmetric $2$-$(35,17,8)$ designs ${\mathcal D}$ with an automorphism of order 2, the number of which is given in the third row, have the full automorphism group of size greater than 2. \\

\begin{table}[htpb!]
\begin{center} \begin{footnotesize}
\begin{tabular}{|c ||c| c | c|c|c|c| }
 \hline 
 
Number of fixed points& 3 & 5 & 7 & 9 & 11 &15 \\
 \hline
Number of  OM & 9007& 3752 & 3291 & 275 & 1506 & 1063  \\
\hline
Number of  designs&111098 & 237058 &687649 & 430654 &8386387 &1817486 \\
\hline 
Designs with $ |Aut({\mathcal D})| > 2$ &5237 & 4301& 16811&3116 & 15945& 10444 \\
\hline
\end{tabular} \end{footnotesize}
 \caption{$2$-$(35,17,8)$  designs with an automorphism of order 2}\label{OM}
\end{center} 
\end{table}

If the full automorphism group of a symmetric $2$-$(35,17,8)$ design is of order 2, then such a design admits only one action of an automorphism of order 2. However, if the full automorphism group is of composite order and contains different conjugacy classes of  subgroups of order 2, it is possible that such a design admits different actions of an automorphism of order 2. We therefore considered all obtained designs with a full automorphism group of size greater than 2, the number of which is given in the forth row of Table \ref{OM}, and performed an isomorphism test. 
Our analysis leads to a complete classification, up to isomorphism, of the symmetric $2$-$(35,17,8)$ designs admitting an automorphism of order 2, and  is summarized in Theorem \ref{des35}.\\

\begin{thm} \label{des35}

\begin{itemize}
\item An automorphism of order 2 of a symmetric $2$-$(35,17,8)$ design fixes $f$ points, where $f \in \{3,5,7,9,11,15 \}$.
\item  Up to isomorphism, there are exactly $11,642,495$ symmetric $2$-$(35,17,8)$ designs admitting an automorphism 
of order 2.
\item The full automorphism groups of the symmetric  $2$-$(35,17,8)$ designs  with an automorphism of order 2
 are described in Table  \ref{tab_inv_all}.  

\end{itemize}

\end{thm}

\begin{table}[htpb!] 
\begin{center} \begin{scriptsize}
\begin{tabular}{|c |c| c ||c|c|c|}
 \hline   
 
The order& The structure&No. of &The order& The structure&No. of \\
of $Aut\mathcal{(D)}$&of $Aut\mathcal{(D)}$&designs&of $Aut\mathcal{(D)}$&of $Aut\mathcal{(D)}$&designs\\
 \hline \hline

 40320& $S_8$&1 & 24& $Z_2\times A_4$ &27\\
\hline
  1152&$((((E_8:Z_3):Z_2):Z_3):Z_2):Z_2$ &1 & & $SL(2,3)$ &4\\
\hline
 720&$S_3\times S_5$ &1 & & $Z_3\times D_8$ &2\\
\hline
 420&$Z_5:(Z_7:Z_{12})$ &1 & & $E_4\times S_3$ &5\\
\hline
 384&$(((E_8:E_4):Z_3):Z_2):Z_2$ &1 & &$S_4$  &7\\
\hline
 288&$(A_4\times A_4):Z_2$ &1 & 18& $Z_3\times S_3$ &26\\
\hline
 192&$((E_8:E_4):Z_3):z_2$ &1 & &$E_9:Z_2$  &15\\
\hline
 168&$E_8:(Z_7:Z_3)$ &1 & &$D_{18}$  &2\\
\hline
 144&$S_3\times S_4$ &1 & 16& $Z_2 \times D_8$  &32\\
\hline
 &$Z_2\times ((S_3\times S_3):Z_2)$ &1 & &$(Z_4\times Z_2) :Z_2$  &10\\
\hline
 136&$Z_{17}:Z_8$ &1 & &$E_{16}$  &60\\
\hline
 128&$(D_8\times D_8):Z_2$ &1 & & $(Z_4\times Z_2):Z_2$ &12\\
\hline
 96&$(E_{16}:Z_2):Z_3)$ &2 & &$Z_4:Z_4$  &2\\
\hline
 &$D_8\times A_4$ &2 & &$QD_{16}$  &1\\
\hline
 72& $Z_2\times(E_9:Z_4)$ &2 & &$Z_4\times E_4$  &2\\
\hline
 &$Z_2\times S_3\times S_3$ &1 & & $D_{16}$ &2\\
\hline
 &$(S_3\times S_3):Z_2$ &1 & 12&$D_{12}$  &47\\
\hline
 &$(Z_3\times A_4):Z_2$ &2 & & $Z_{12}$ &14\\
\hline
 64&$(E_{16}:Z_2):Z_2$ &1 & & $Z_6\times Z_2$ &62\\
\hline
 &$D_8\times D_8$ &2 & & $A_4$ &29\\
\hline
 &$(E_8:Z_4):Z_2$ &2 & & $Z_3:Z_4$ &4\\
\hline
 60&$Z_3\times(Z_5:Z_4)$ &1 & 10& $Z_{10}$  &2\\
\hline
 48&$Z_2\times S_4$ &3 & &$D_{10}$  &1\\
\hline
 &$GL(2,3)$ &1 & 8& $D_8$ &333\\
\hline
 40& $Z_2\times (Z_5:Z_4)$ &1 & & $E_8$ &493\\
\hline
 36& $S_3\times S_3$&3 & &$Z_8$  &139\\
\hline
 &$E_9:Z_4$ & 6& &$Z_4 \times Z_2$  &202\\
\hline
 & $Z_2\times (E_9:Z_2)$&2 & &$Q_8$  &58\\
\hline
 &$Z_6\times S_3$ &2 & 6& $S_3$ &249\\
\hline
 32&$E_{16}:Z_2$ &12 & & $Z_6$ &1667\\
\hline
 &$E_8:Z_4$ &5 & 4& $E_4$  &20171\\
\hline
 &$E_4\times D_8$ &4 & & $Z_4$ &4266\\
\hline
 24& $(Z_6\times Z_2):Z_2$&4 & 2&$Z_2$  &11614478\\
\hline

\end{tabular} \end{scriptsize} 
 \caption{The full automorphism groups of the symmetric $2$-$(35,17,8)$ designs with an automorphism of order 2} \label{tab_inv_all}
\end{center} 
\end{table}

\section{Related $3$-$(36,18,8)$ designs} \label{3des}

Every $2$-$(4t-1,2t-1,t-1)$ design $D$ is extendable in a unique way (up to isomorphism) to a $3$-$(4t,2t,t-1)$ design $D^*$, called Hadamard $3$-design, by adding one new point to all blocks
of $D$, and adding $4t-1$ new blocks of size $2t$ being the complements of the blocks of $D$. This process can be reversed with respect to every of the $4t$ points of $D^*$ to produce a "derived" $2$-$(4t-1,2t-1,t-1)$ design. Two derived $2$-$(4t-1,2t-1.t-1)$ designs of a $3$-$(4t,2t,t-1)$ design $D^*$ obtained from $D^*$ with respect to points $i$ and $j$ are isomorphic if and only if there is an automorphism of $D^*$ which maps point $i$ to point $j$.\\

In this section we discuss the $3$-$(36,18,8)$ designs obtained from the symmetric $2$-$(35,17,8)$ designs with an automorphism of order 2.
After constructing all such designs, we eliminate isomorphic copies 
and obtain $1,015,225$ pairwise nonisomorphic  $3$-$(36,18,8)$ designs.
Information about the full automorphism groups of these designs is given in Table \ref{tab_des3_all}.

\begin{table}[htpb!] 
\begin{center} \begin{scriptsize}
\begin{tabular}{|c |c| c ||c|c|c|}
 \hline   
 
The order& The structure&No. of &The order& The structure&No. of \\
of $Aut\mathcal{(D)}$&of $Aut\mathcal{(D)}$&designs&of $Aut\mathcal{(D)}$&of $Aut\mathcal{(D)}$&designs\\
 \hline \hline

 40320& $S_8$&1 & 24& $Z_2\times A_4$ &14\\
\hline
 720&$S_3\times S_5$ &1 & & $SL(2,3)$ &6\\
\hline
 576&$(A_4 \times A_4) : Z_4$ &1 & & $Z_3\times D_8$ &1\\
\hline
432&$((E_9:Q_8):Z_3):Z_2$&1& & $(Z_6\times Z_2):Z_2$&2\\
\hline
 420&$Z_5:(Z_7:Z_{12})$ &1 & & $E_4\times S_3$ &1\\
\hline
 384&$(((E_8:E_4):Z_3):Z_2):Z_2$ &1 & &$S_4$  &4\\
\hline
 272&$Z_{17}:Z_{16}$  &1 & 18& $Z_3\times S_3$ &28\\
\hline
 192&$((E_8:E_4):Z_3):Z_2$ &1 & &$E_9:Z_2$  &3\\
\hline
 168&$E_8:(Z_7:Z_3)$ &1 & 16&$Z_{16}$  &10\\
\hline
 144& $E_9:QD_{16}$ &1 & & $Z_2 \times D_8$  &6\\
\hline
 &$Z_2\times ((S_3\times S_3):Z_2)$ &1 & &$(Z_4\times Z_2) :Z_2$  &5\\
\hline
 108&$(E_9 : Z_3):E_4$ &1 & &$E_{16}$  &30\\
\hline
 96 &$(Z_4 \times Z_4 :Z_2) :Z_3$  &3 & & $(Z_4\times Z_2):Z_2$ &5\\
\hline
 &$(E_{16}:Z_2):Z_3$ &1 & &$Z_4:Z_4$  &1\\
\hline
 &$D_8\times A_4$ &1 & &$QD_{16}$  &18\\
\hline
 72& $Z_2\times(E_9:Z_4)$ &2 & &$Z_4\times E_4$  &1\\
\hline
 &$Z_2\times S_3\times S_3$ &1 & & $D_{16}$ &1\\
\hline
 &$(Z_3\times A_4):Z_2$ &1& &$Q_{16}$&3 \\
\hline
 64&$( (Z_8:Z_2):Z_2):Z_2$ &1 & 12&$D_{12}$  &26\\
\hline
 &$D_8\times D_8$ &1 & & $Z_{12}$ &6\\
\hline
 &$(E_8:Z_4):Z_2$ &3 & & $Z_6\times Z_2$ &44 \\
\hline
 60&$Z_3\times(Z_5:Z_4)$ &1 & & $A_4$ &25\\
\hline
 54&$Z_9:Z_6$ &2 & & $Z_3:Z_4$ &2\\
\hline
&$(E_9 : Z_3) :Z_2$&2& 10& $Z_{10}$  &2\\
\hline
 48&$GL(2,3)$ &3 & &$D_{10}$  &1 \\
\hline
 40& $Z_2\times (Z_5:Z_4)$ &1 & 8& $D_8$ &102\\
\hline
 36&$Z_6\times S_3$ &2 & & $E_8$ &36\\
\hline
 &$E_9:Z_4$ &3 & &$Z_8$  &49\\
\hline
 & $Z_2\times (E_9:Z_2)$&2 & &$Z_4 \times Z_2$  &70\\
\hline
 32&$Z_{16}:Z_2$  &1 & &$Q_8$  &50\\
\hline
 &$E_{16}:Z_2$ &4 & 6& $S_3$ &82\\
\hline
 &$E_4.E_8$ &2 & & $Z_6$ &696\\
\hline
 &$E_4\times D_8$ &2 & 4& $E_4$  &3792\\
\hline
 &$(Z_4 \times Z_4):Z_2$ &2 & & $Z_4$ &1375\\
\hline
&$QD_{32}$&1& 2&$Z_2$  &1008672 \\
\hline
&$(Z_8 \times Z_2) : Z_2$&1& & &\\
\hline
&$Q_{32}$&1& & &\\
\hline

\end{tabular} \end{scriptsize} 
 \caption{The full automorphism groups of the obtained $3$-$(36,18,8)$ designs} \label{tab_des3_all}
\end{center} 
\end{table}
 
\bigskip
\noindent {\bf Acknowledgement} 

The first author is supported by the University of Rijeka project uniri-iskusni-prirod-23-62.\\



\end{document}